\documentclass[fleqn,12pt]{article}
\usepackage{amssymb,amsmath,amsthm}
\usepackage{amsthm}
\usepackage{a4}
\usepackage{amssymb}
\usepackage{amsmath}

\newtheorem{theorem}{\bf Theorem}[section]
\newtheorem{lemma}{\bf Lemma}[section]

\newcommand{\vol}{\mbox{vol}\,}

\newcommand{\C}{{\mathbb C}}

\newcommand{\Gm}{{\mathbb G}_{\rm m}}

\newcommand{\Z}{{\mathbb Z}}

\newcommand{\R}{{\mathbb R}}

\newcommand{\bea}{\begin{eqnarray*}}
\newcommand{\eea}{\end{eqnarray*}}
\newcommand{\be}{\begin{eqnarray}}
\newcommand{\ee}{\end{eqnarray}}

\newcommand{\ve}{\boldsymbol}

\newcommand{\supp}{\mbox{supp}}
\newcommand{\spn}{{\rm span}}

\newcommand{\power}{z}

\begin{document}

\setlength{\unitlength}{1mm}
\setcounter{page}{1}

\title{Power maps and subvarieties of the complex algebraic $n$--torus
}
\author{Iskander Aliev and Chris Smyth}

\maketitle

\begin{abstract}

Given a subvariety $V$ of the complex algebraic torus $\Gm^n$ defined by polynomials of total degree at most $d$ and a power map $\phi: \Gm^n \rightarrow \Gm^n$, the points ${\ve x}$ whose forward orbits ${\mathcal O}_\phi({\ve x})$
belong to $V$ form its {\em stable} subvariety $S(V,\phi)$. The main result of the paper provides an upper bound $T=T(n,d,\phi)$ for the number of iterations of the power map $\phi$ required to ``cut off'' the points of $V$ that do not belong to $S$.

\bigskip
\noindent {\bf 2000 MS Classification:} Primary 11G35; Secondary
14L40.
\end{abstract}

\section{Introduction}

Given a set $M$ and a map $\phi: M \rightarrow M$, one of the goals of dynamics is to classify the points ${\ve x}$ of $M$
according to the behavior of their forward orbits  ${\mathcal O}_\phi({\ve x})$.
Let $\Gm^n$ be the complex algebraic $n$-torus and
$\power$ be an integer $\ge 2$.
This paper is concerned with a special case of the general classification problem for $M=\Gm^n$ and the {\em
$\power$-th power map} $\phi\,: \Gm^n\rightarrow\Gm^n$ defined by
the rule
\bea
\phi\,: \; (x_1,\ldots,x_n)\longmapsto (x_1^\power,\ldots,x_n^\power)\,.
\eea
We refer the reader to Silverman \cite{Silverman} and Zhang \cite{Zhang} for results and further references on polynomial maps in algebraic dynamics.

As an affine variety, we will identify the torus $\Gm^n$ with the
Zariski open subset $x_1x_2\cdots x_n\neq 0$ of affine space
${\mathbb A}^n$, with the usual multiplication
\bea (x_1,x_2,\ldots,x_n) \cdot
(y_1,y_2,\ldots,y_n)=(x_1y_1,x_2y_2,\ldots,x_ny_n)\,. \eea
By {\em algebraic subvariety}, or simply {\em subvariety} of
$\Gm^n$, we will understand a Zariski closed subset. An {\em
algebraic subgroup} of $\Gm^n$ is a Zariski closed subgroup. A
{\em subtorus} of $\Gm^n$ is a geometrically irreducible algebraic
subgroup.
Recall that the torsion points of
$\Gm^n$ are precisely the points ${\ve \omega}=(\omega_1, \ldots,
\omega_n)$ whose coordinates $\omega_i$ are roots of
unity.
By {\em torsion coset} we will mean a coset ${\ve \omega}H$, where
$H$ is a subtorus of $\Gm^n$ and ${\ve \omega}$ is a torsion point.
Both for convenience and for technical reasons we will often use Laurent
polynomials through this paper. For a Laurent polynomial
\bea f({\ve X})=\sum_{{\ve i}=(i_1,\ldots,i_n)\in \Z^n}a_{\ve
i}{\ve X}^{\ve i}\,,\;\;{\ve X}^{\ve i}=X_1^{i_1}\cdots
X_n^{i_n}\,,\eea
we will write
 \bea f^t=f(\phi^t({\ve X}))=f(X_1^{\power^t},\ldots,X_n^{\power^t})\,.\eea
The notation $f \sim g$ for nonzero Laurent polynomials $f$
and $g$ will mean that $f/g$ is either a constant or a
monomial. A Laurent polynomial $f$ will be called {\em nontrivial}
if $f\nsim 1$. The same notation will be applied to standard
polynomials.

Now let $f_1, f_2, \ldots, f_s\in \C[X_1,\ldots,X_n]$ be
polynomials of total degree $\le d$. For the subvariety $V=Z(f_1,
f_2, \ldots, f_s)\subset\Gm^n$ and an integer $u \ge 0$ we define
the set
\bea V(\power, u)=\{{\ve x}\in \Gm^n: \phi^t({\ve x})\in
V\;\;\mbox{for}\;\;t=0,\ldots,u\}\,. \eea
Clearly,  $V(\power, u)$ is the subvariety of $\Gm^n$ defined by
the polynomials $f_i^t$ for $i=1,\ldots, s$ and $t=0,\ldots,u$.
A subvariety $S\subset V$ will be called ($\power, V$)-{\em
stable} if for all positive integers $u$ we have $S\subset
V(\power, u)$. In other words, for any point ${\ve x}$ of the
subvariety $S$ its forward orbit ${\mathcal O}_\phi({\ve x})$
belongs to $V$. We will call a $(\power,
V)$-stable subvariety $S$ {\em maximal} if there is no $(\power,
V)$-stable subvariety $S'$ with $S\varsubsetneq S'$.
 The maximal $(\power, V)$-stable subvariety will be denoted by $S(\power, V)$.

If for some integer $T$ the subvariety $S=V(\power, T)$ is $(\power,
V)$-stable then clearly $S=S(\power, V)$.
The main result of the paper states that such an integer $T$ can be
effectively  bounded in terms of the power $\power$, the dimension $n$
and the maximum total degree $d$ of the defining polynomials $f_1,
f_2, \ldots, f_s$ of the subvariety $V$.

%
%
\begin{theorem}
There are effectively computable constants $T=T(\power, n,d)$,
$E=E(\power, n,\,d)$ and $L=L(\power, n,\,d)$ such that

\begin{itemize}

\item[(i)] for any  subvariety $V$ of $\Gm^n$ defined by the
polynomials of total degree at most $d$ we have
$S(\power, V)=V(\power, T)$;

\item[(ii)] the subvariety $S(\power, V)$ is
contained in a finite union of $(n-1)$-dimensional torsion cosets
$\bigcup_i D_i$, where
\bea D_i=Z(h)\,,\;\;h\sim {\ve X}^{{\ve a}_i}-\zeta_i\eea
with $||{\ve a}_i||_2\le L$, $\zeta_i$  a $\power^k(\power^l-1)$th
root of unity, and $k,l\le E$.

\end{itemize} \label{Intersections}
\end{theorem}

For the sake of completeness in Section \ref{main_proof} we will give recurrent formulae for
all the constants involved in the main theorem.

We will now discuss a relation between Theorem \ref{Intersections} and the classical theorem of Skolem--Mahler--Lech.
Let us consider the sequence $\{u_m\}$ defined as
\bea
u_m=a_1 \alpha_1^m+\cdots+ a_N \alpha_N^m\;\;\; (m\in\Z)
\eea
with nonzero coefficients $a_i\in \C$ and with nonzero distinct elements $\alpha_i\in \C$. The sequence $\{u_m\}$ is called {\em nondegenerate}
if no quatient $\alpha_i/\alpha_j$ $(1\le i<j\le N)$ is a root of unity.
Let $S(u_m)$ denote the set of zeros of $\{u_m\}$, i.e., the set of solutions $k\in \Z$ of the equation  $u_k=0$.
The Skolem--Mahler--Lech theorem implies that for nondegenerate sequences $\{u_m\}$ the set $S(u_m)$ is finite.
Furthermore, Theorem 1.2 of Evertse, Schlickewei and Schmidt \cite{ESS} gives an upper bound for the cardinality of $S(u_m)$ in terms of
$N$ only.

Let $f({\ve X})=\sum_{{\ve i}}a_{\ve
i}{\ve X}^{\ve i}$ be
a polynomial of total degree $d$ and $H=Z(f)\subset\Gm^n$ be the hypersurface defined by $f$.
Every point ${\ve \beta}$ of the stable subvariety $S(\power, H)$ clearly corresponds to a sequence $\{u_m\}$ with infinite
set $S(u_m)$. Moreover, the number of terms $N$ of the sequence $\{u_m\}$ will depend on $n$ and $d$ only. Consequenly, the above mentioned theorem of Evertse, Schlickewei and Schmidt \cite{ESS} gives us a restriction on the points that can lie on $S(\power, H)$. Indeed, it implies that the ratio of at least two of the numbers
${\ve \beta}^{\ve i}$ (${\ve i}$ are integer vectors such that $a_{\ve i}\neq 0$) must be a root of unity. Moreover, the upper bound for the cardinality of $S(u_m)$ gives an upper bound for the number of iterations required to cut off the points that do not satisfy this restriction.

The next theorem shows that for $n\le 2$ the maximal stable subvarieties have remarkably simple structure.


\begin{theorem}
Let $V$ be a subvariety of $\Gm^n$, $n\le 2$,  and $\power$ be an integer $\ge 2$.
Then the maximal $(\power, V)$-stable subvariety $S(\power, V)$ is a finite union of
torsion cosets.
\label{2D}
\end{theorem}

The general problem of finding all torsion cosets on a given subvariety of $\Gm^n$ has been addressed in Aliev and Smyth \cite{Aliev-Smyth}, Ruppert \cite{Ruppert} and Sarnak and Adams \cite{Sarnak}. For the special case $n=2$ the best among known algorithms is due to Beukers and Smyth \cite{Beukers-Smyth}.
Combining Theorem \ref{2D} and the algorithm of Beukers and Smyth we can find maximal stable subvarieties for the plane curves. For instance, let us consider the following example.

\vskip.3cm
\noindent {\bf Example}: Let $f(X,Y)=X^2Y^2+X^2Y+XY^2+XY+X+Y+1$. The curve $V=Z(f)$ has the largest known ratio
\bea
\frac{N_{tor}(f)}{\vol(f)}=16\,,
\eea
where $N_{tor}(f)$ denote the number of torsion points on the curve $Z(f)$ and $\vol(f)$ is the volume of the Newton polygon of the polynomial $f$ (see Beukers and Smyth \cite{Beukers-Smyth} for details). Indeed $Z(f)$ contains exactly $48$ torsion points of $\Gm^2$ listed below:
\begin{itemize}
\item[] $(\omega_{12}^{4i},\omega_{12}^i)$, $(\omega_{12}^{i},\omega_{12}^{4i})$, $(-\omega_{12}^{i},\omega_{12}^{i})$,\;\; $i=1,5,7,11$\,;
%
%
\item[] $(\omega_{7}^{i},\omega_{7}^{2i})$, $(\omega_{7}^{2i},\omega_{7}^{i})$,\;\; $i=1,\ldots,6$\,;
\item[] $(-\omega_{30}^{3i},\omega_{30}^i)$, $(\omega_{30}^{i},-\omega_{30}^{3i})$, $(\omega_{30}^{i},\omega_{30}^{11i})$,\;\; $i=1,7,11,13,17,19,23,29$\,.
\end{itemize}
The power representation of torsion points lying on $V$ allows us to easily detect the maximal $(\power, V)$-stable subvarieties. For instance,
it is clear that $S(2, V)$ consists of $12$ torsion points $(\omega_{7}^{i},\omega_{7}^{2i})$, $(\omega_{7}^{2i},\omega_{7}^{i})$,\;\; $i=1,\ldots,6$.

\section{Lattices, algebraic subgroups of $\Gm^n$ and geometry of numbers}

%
We recall some basic definitions. A {\em lattice} is a discrete
subgroup of $\R^n$. Given a lattice $L$ of rank $r$, any set of
vectors $\{{\ve b}_1,\ldots,{\ve b}_r\}$ with $L={\rm
span}_{\Z}\{{\ve b}_1,\ldots,{\ve b}_r\}$ or the matrix
${\boldsymbol B}=({\ve b}_1,\ldots,{\ve b}_r)$ with rows ${\ve
b}_i$ will be called a {\em basis} of $L$.  The {\em determinant}
of a lattice $L$ with a basis ${\boldsymbol B}$ is defined to be
\bea  \det (L) = \sqrt{\det({\boldsymbol B}\,{\boldsymbol B}^T)}\,.
\eea
%
%
%
When $L$ is a lattice of rank $n$, its {\em polar} lattice $L^*$
is defined as
\bea L^*=\{{\ve x}\in \R^n : \langle {\ve x},{\ve
y}\rangle\in\Z\;\;\mbox{for all}\;\; {\ve y}\in L\}\,, \eea
where $\langle \cdot,\cdot\rangle$ denotes the inner product.
 Given a basis ${\mathcal B}=({\ve b}_1,\ldots,{\ve b}_n)$ of
$L$, the basis of $L^*$ {\em polar} to ${\mathcal B}$ is the basis
${\mathcal B}^*=({\ve b}_1^*,\ldots,{\ve b}_n^*)$ such that
\bea \langle {\ve b}_i\,, {\ve b}_j^*\rangle=\delta_{ij}\,,\;\;
i,j=1,\ldots, n\,,\eea
with $\delta_{ij}$ the Kronecker delta.

By an {\em integer} lattice we understand a lattice
$A\subset\Z^n$. An integer lattice is called {\em primitive} if
$A=\spn_{\R}(A)\cap \Z^n$. Here $\spn_{\R}(A)$ denotes the
subspace of $\R^n$ spanned by the vectors of the lattice $A$. For
an integer lattice $A$, we define the subgroup $H_A$ of $\Gm^n$ by
\bea H_A=\{{\ve x}\in\Gm^n: {\ve x}^{\ve a}=1\;\; {\rm for\;
all}\;\; {\ve a}\in A\}\,. \eea
Then, for instance, $H_{\Z^n}$ is the trivial subgroup.
\begin{lemma}
The map $A\mapsto H_A$ sets up a bijection between integer
lattices and algebraic subgroups of $\Gm^n$. A subgroup $H=H_A$ is
irreducible if and only if the lattice $A$ is primitive.
\label{lattices_cosets}
\end{lemma}
\begin{proof}
See Lemmas 1, 2 of Schmidt \cite{Schmidt}.
\end{proof}

Let now ${\boldsymbol B}=({\ve b}_1, {\ve b}_2,\ldots, {\ve b}_n)$ be
a basis of the lattice $\Z^n$. The map $\psi: \Gm^n\rightarrow
\Gm^n$ defined by
\be \psi({\ve x})=({\ve x}^{{\ve b}_1},\ldots, {\ve x}^{{\ve
b}_{n}})\, \label{auto} \ee
is an automorphism of $\Gm^n$ (see Ch. 3 in Bombieri and Gubler
\cite{Bombieri-Gubler} and Section 2 in Schmidt \cite{Schmidt}).
The automorphism (\ref{auto}) is traditionally called a {\em
monoidal transformation}.  To make the inductive argument used in
the proof of Theorem \ref{Intersections} more transparent, we will
associate with ${\boldsymbol B}$ the new coordinates
$(Y_1,\ldots,Y_n)$ in $\Gm^n$ defined by
\be Y_1={\ve X}^{{\ve b}_1},\;\; Y_2={\ve X}^{{\ve b}_2},
\ldots,\;\; Y_{n}={\ve X}^{{\ve b}_{n}}\,.
\label{associated_coordinates}\ee
Thus changing variables ${\ve X}\mapsto{\ve Y}$ is equivalent to applying the automorphism $\psi$ defined by (\ref{auto}).

Suppose that the matrix ${\boldsymbol B}^{-1}$ has rows ${\ve
r}_1, {\ve r}_2,\ldots, {\ve r}_n$. By the {\em image} of a
Laurent polynomial $f\in \C[X_1^{\pm 1},\ldots,X_n^{\pm 1}]$ in
coordinates $(Y_1,\ldots,Y_n)$ we mean the Laurent polynomial
\bea f^{\boldsymbol B}({\ve Y})=f({\ve Y}^{{\ve r}_1},\ldots,{\ve
Y}^{{\ve r}_n})\,. \eea
and by the {\em image} of a subvariety $V=Z(f_1,\ldots,f_s)$ we
understand the subvariety $V^{\boldsymbol B}$ defined as the set of common zeroes in $\Gm^n$ of the Laurent polynomials $f_1^{\boldsymbol
B},\ldots,f_s^{\boldsymbol B}$.
Then the subvariety $V^{\boldsymbol B}$ is well defined and we will simply write $V^{\boldsymbol B}=Z(f_1^{\boldsymbol
B},\ldots,f_s^{\boldsymbol B})$.

Let $C={\ve \omega}H_A$ be an $(n-r)$-dimensional torsion coset and let $({\ve a}_1, {\ve a}_2,\ldots, {\ve a}_{r})$ be a
basis of the lattice $A$. Then $C$ can be defined by $r$ equations
\bea {\ve X}^{{\ve a}_i}-{\ve \omega}^{{\ve a}_i}=0\,,\;\;\;
i=1,\ldots,r\,.\eea
Each such an equation defines an $(n-1)$-dimensional torsion coset in $\Gm^n$ and, after a suitable automorphism, will have the form
$X_i=\omega_i$ with $\omega_i$ a root of unity. This observation allows us to reduce the dimension of the problem provided that
the subvarieties of interest lie in a torsion coset.

Let $G^\bot$ denote the orthogonal complement of the subspace
$G\subset\R^n$. We will need the following technical lemma from
geometry of numbers.
\begin{lemma}
Let $G$ be a subspace of $\R^n$ with $\dim(G)={\rm
rank}(G\cap\Z^n)=r<n$. Then there exists a basis ${\boldsymbol
A}=({\ve a}_1, {\ve a}_2,\ldots,{\ve a}_n)$ of the lattice $\Z^n$
such that ${\ve a}_1\in G^\bot$ and the vectors of the polar basis
${\boldsymbol A}^*=({\ve a}^*_1, {\ve a}^*_2,\ldots,{\ve a}^*_n)$
satisfy the inequalities
\be
||{\ve a}^*_i||_2< 1+\frac{n-1}{2}\gamma_{n-1}^{\frac{n-1}{2}}
\gamma_{n-r}^{\frac{1}{2}}\det(G\cap\Z^n)^{\frac{1}{n-r}}\,,\;\;\;
i=1,\ldots,n\,.\ee
with
$\gamma_{n-1}$  the Hermite constant for dimension $n-1$ (for the definition and properties of the Hermite constant see
Section 38.1 of Gruber--Lekkerkerker \cite{Gruber-Lekkerkerker}).
\label{upper_bound_for_inverse_elements}
\label{suitable_matrix}
\end{lemma}

\begin{proof}
See Corollary 2.1 of Aliev and Smyth \cite{Aliev-Smyth}.
\end{proof}

\section{Proof of Theorem \ref{Intersections}}
\label{main_proof}

In this section we will show that the statement of Theorem \ref{Intersections} holds for the constants $T$, $E$ and $L$ defined as follows.
We can take $T(\power,1,d)=d$ and
%
%
\bea T(\power,n,d)= c_0(\power, d) T(\power^{c_0(\power, d)}, n-1, c_1(\power,
n,d)\power^{2 c_0(\power, d)-1} )\eea \bea+c_0(\power, d)\,,\; \;\eea where $c_0(\power, d)=\lceil3\log(d)/\log(\power)\rceil+d(d-1)$ and $c_1(\power,
n,d)=n(n+1)d+\lceil n(n^2-1)\gamma_{n-1}^{\frac{n-1}{2}} d L(\power, n, d)\rceil$.
The constants $E$ and $L$ can be calculated as
\bea E(\power, n,\,d)= E(\power, n-1, d^5\power^{d(d-1)})\,,\; E(\power, 1,d)=d\,;\eea
\bea L(\power, n,\,d)= L(\power, n-1, d^5 \power^{d(d-1)})\,, \;n\ge 3\,,\eea\bea
L(\power, 2,\,d)= 2d\,,\; L(\power, 1,\,d)=1\,. \eea

It is clearly enough to prove the theorem in the case when $V$ is
a hypersurface in $\Gm^n$.
We will start with the following auxiliary result which slightly
extends the Proposition 1 of Ruppert \cite{Ruppert}.

\begin{lemma}
Let $f\in \C[X_1, X_2,\ldots, X_n]$, $n\ge 2$, be a nontrivial irreducible
polynomial and $\eta_1, \ldots, \eta_n$ be the $M$th roots of
unity. If $f$ divides the polynomial $f(\eta_1 X_1^N, \eta_2 X_2^N, \ldots,
\eta_n X_n^N)$ for some integer $N>1$ then $f$ has the form
\bea f \sim {\ve X}^{\ve a}- \zeta \,,\eea
where ${\ve a}$ is an integer vector and $\zeta$ is a root of
unity with $\zeta^{(N-1)M}=1$.

\label{Ruppert-plus}
\end{lemma}

\begin{proof}
We will modify the proof of the Proposition 1 of Ruppert
\cite{Ruppert}. It consists of seven steps and only the first and
the last of them need to be changed.

The first step allows us to assume that the polynomial $f$ cannot
be represented in the form $f=g(X_1^{a_1},\ldots,X_n^{a_n})$ with
some $a_i>1$. To see this observe that in our case the polynomial
$g$ will divide the polynomial $g(\eta_1^{a_1} X_1^N, \eta_2^{a_2}
X_2^N, \ldots, \eta_n^{a_n} X_n^N)$ and $\eta_1^{a_1},
\eta_2^{a_2}, \ldots, \eta_n^{a_n}$ are the $M'$th roots of unity
with $M'|M$. Therefore we can simply replace $f$ by $g$.

The next five steps of the original proof will remain almost
unchanged - we just always replace the polynomial $f(X_1^N, X_2^N,
\ldots, X_n^N)$ by the polynomial $f(\eta_1 X_1^N, \eta_2 X_2^N,
\ldots, \eta_n X_n^N)$.

At the seventh step we may assume that $f=X_1\cdots X_m- d
X_{m+1}\cdots X_n$ and, consequently, $f(\eta_1 X_1^N, \eta_2
X_2^N, \ldots, \eta_n X_n^N)= \mu_1 (X_1\cdots X_m)^N - d
\mu_2(X_{m+1}\cdots X_n)^N$, where $\mu_i$ are the $M$th roots of
unity. After the substitution $X_2=\mu_1^{-1/N}, X_3=1,\ldots,
X_{n-1}=1, X_n=\mu_2^{-1/N}$, we have $\mu_1^{-1/N}X_1- d
\mu_2^{-1/N} | X_1^N-d$. Therefore, $d^{(N-1)M}=\mu_2^M=1$.

\end{proof}

The next lemma show that, roughly speaking, for all sufficiently large $t$ the irreducible common factors of the polynomials $f, f^1, f^2, \ldots, f^t$
will define the $(n-1)$-dimensional torsion cosets.

\begin{lemma}

Let $f\in \C[X_1,\ldots,X_n]$, $n\ge 2$, be a (possibly reducible) polynomial of the total degree $d$. Then
\begin{itemize}
\item [(i)] for any integer $t\ge
t_0=\lceil3\log(d)/\log(\power)\rceil$ each irreducible common
factor of the polynomials $f$ and  $f^t$ consists of precisely two
terms;
\item [(ii)] there exists an integer, $t$, with $t_0\le t_1\le c_0=c_0(\power, d)$
such that each irreducible common factor $g$ of the polynomials
$f$ and $f^{t_1}$ has the form
\be g \sim {\ve X}^{\ve a}- \zeta
\,,\label{polynomial-of-coset}\ee
where ${\ve a}$ is an integer vector and $\zeta$ is a
$\power^k(\power^{l}-1)$th root of unity with $0\le k,l\le
c_0$.

\end{itemize}

\label{Gcd}
\end{lemma}

\begin{proof}

(i) Suppose that the statement is wrong. Then for some irreducible
factors $g_1$ and $g_2$ of the polynomial $f$ we would have
$g_1|g_2^t$ with $t\ge t_0$ and $g_1$ consisting of more than two terms.
Furthermore, by the first result stated in Section III
of Gourin \cite{Gourin}, the polynomial $g_2$ must consist of more than two terms as well.

Observe that $\deg(g_2^t)\ge \power^{t_0}\ge d^3$ and thus $g_1\neq g_2^t$, so that the polynomial $g_2^t$ is reducible.
By Theorem I of Gourin \cite{Gourin}
there exist nonnegative integers $t_1, t_2,\ldots,t_n$ such that
\be
\power^{t_i}\le d^2\,, \;\;i=1,\ldots, n\,
\label{bound-for-t_i}
\ee
and
$g_1$ can be obtained by replacing each $X_i$ by $X_i^{\power^{t-t_i}}$ in an irreducible factor $g_2^*$ of the polynomial $g_2(X_1^{\power^{t_1}},\ldots, X_n^{\power^{t_n}})$.
Therefore we must have $\deg(g_1)> \power^{t-t_i}$. However by (\ref{bound-for-t_i})
\bea
\power^{t-t_i}\ge \frac{\power^{t_0}}{d^2}\ge d\,.
\eea
The  contradiction obtained proves part (i).

(ii) The polynomial $f$ can be written in the form
\bea f=h_1 h_2 \cdots h_k f'\,,\eea
where each of  the polynomials $h_i$ consists of precisely two terms and
each of  factors of the polynomial $f'$ have more than two terms. Let
us choose an integer $s_0$ with $t_0\le s_0\le c_0$. By part (i) and the first result stated in Section III
of Gourin \cite{Gourin}, for
some $1\le i,j\le k$ we have $h_i|h_j^{s_0}$. Put $g=h_i$ and suppose
first that $i=j$ so that $g|g^{s_0}$. Then by Proposition 1 of
Ruppert \cite{Ruppert} the polynomial $g$ has the form
(\ref{polynomial-of-coset}) with $\zeta^{\power^{s_0}-1}=1$. Thus we can take $t_1=l=s_0$ and $k=0$. To settle
the case $i\neq j$ suppose that $g|h_j^{s_1}$ for some integer
$s_1>s_0$. We will show that this implies that the polynomial $g$
has the desired form. By Lemma I of Gourin \cite{Gourin} we have
\bea h_j^{s_1}= \prod_v g_v^{s_1-s_0}\,, \eea
where the polynomials $g_v$ are irreducible and form the {\em
complete set of $\power^{s_0}\cdots \power^{s_0}$ transforms}, as
introduced in Gourin \cite{Gourin} p. 486, obtained from the
polynomial $g$.
Therefore  $g$ divides a polynomial
$g(\eta_1 X_1^{\power^{s_1-s_0}}, \eta_2 X_2^{\power^{s_1-s_0}}, \ldots,
\eta_n X_n^{\power^{s_1-s_0}})$
where $\eta_i$ are $\power^{s_0}$th roots of unity. Now by Lemma
\ref{Ruppert-plus} the polynomial $g$ must have the form
\bea g \sim {\ve X}^{\ve a}- \zeta \,,\eea
where ${\ve a}$ is an integer vector and $\zeta$ is a root of
unity with $\zeta^{\power^{s_0}(\power^{s_1-s_0}-1)}=1$.

The number $k$ of the factors $h_i$ does not exceed $d$.
Therefore, among the integers $t_0, t_0+1, \ldots, t_0+d(d-1)=c_0$
there always exists an integer $t_1$ satisfying the conditions of
the lemma. Finally we take $k=s_0$ and $l=s_1-s_0$.

\end{proof}

Let now $f\in \C[X_1,\ldots,X_n]$ be a polynomial of total degree
$d$ and $V=Z(f)$. We will proceed by induction on the number of
variables $n$.

At the basis step $n=1$ the defining polynomial $f=f(X)$ has
degree $d$ and $T(\power, 1,d)=d$ as well. Therefore, a root $\omega$ of
$f$ belongs to $S=V(\power, T(\power, 1,d))$ if and only if for some
integers $k\ge 0$, $l>0$ with $k+l\le d$ we have
$\omega^{\power^{k+l}}=\omega^{\power^k}$. Therefore $S$ is
$(\power, V)$-stable and consists of $\power^k(\power^l-1)$th
roots of unity with $k,l\le d=E(\power, 1,d)$. Each such a root $\omega$ can
be regarded as a $0$-dimensional torsion coset $Z(X-\omega)$ in $\Gm$. Thus
$L(\power,1,d)=1$.

Suppose now that $n\ge 2$.
At the inductive step we will first decompose the problem into two
cases. Let us consider the polynomial $g=\gcd(f, f^{t_1})$, where
$t_1$ is the integer satisfying conditions of the part (ii) of
Lemma \ref{Gcd}. The polynomial $g$ can be written as $g=h_1 h_2
\cdots $, where each of the factors $h_i$ is irreducible and has the form
(\ref{polynomial-of-coset}).
Observe that any irreducible component $U$ of the subvariety $S=V(\power, T(\power, n,d))$
belongs to at least one of the subvarieties $Z(g)$ and $Z(f/g,
f^{t_1}/g)$. Thus we can separately consider the case $U\subset
Z(g)$ and the case $U\subset Z(f/g, f^{t_1}/g)\setminus Z(g)$.

Suppose that $U\subset Z(g)$. Since $U$ is irreducible, there is a factor $h_i$ of $g$ with $U\subset Z(h_i)$ and
\bea h_i \sim h={\ve X}^{\ve a}- \zeta\,. \eea
Here $\zeta$ is a $\power^k(\power^{l}-1)$th root of unity with
$0\le k,l\le t_0+d(d-1)$, so that $k, l\le E(\power, n,d)$. Observe
that $||{\ve a}||_1\le 2\deg(h_i)\le 2d$ and thus
\be ||{\ve a}||_2\le 2d \le L(\power, n,d)\label{length-of-a}\,. \ee
Thus the component $U$ lies in the torsion coset $Z(h)$ satisfying the conditions of Theorem \ref{Intersections}.
Let us show that $U$ is $(\power,
V)$-stable.

By Lemma \ref{suitable_matrix} applied to the subspace
$G=(\spn_{\R}({\ve a}))^\bot$, there exists a basis ${\boldsymbol
A}=({\ve a}_1, \ldots, {\ve a}_{n-1})$ of the lattice $\Z^{n}$
such that ${\ve a}_1={\ve a}$ and its polar basis ${\boldsymbol
A}^*=({\ve a}^*_1, \ldots, {\ve a}^*_{n})$ satisfies the
inequality (\ref{upper_bound_for_inverse_elements}).

Let $(Y_1,\ldots,Y_n)$ be the coordinates associated with
${\boldsymbol A}$. Clearly, the subvariety $S$ is $(\power,
V)$-stable if and only if $S^{\boldsymbol A}$ is $(\power,
V^{\boldsymbol A})$-stable.

Consider the subvariety  $X=S^{\boldsymbol A}\cap Z(Y_1-\zeta)$.
The preimage of $X$ in the initial coordinates is the subvariety
$S\cap Z(h)$. Thus $U^{\boldsymbol A}$ is an irreducible component
of $X$ and it  is enough to show that $X$ is $(\power,
V^{\boldsymbol A})$-stable.

Recall that $\zeta$ is a $\power^k(\power^{l}-1)$th root of
unity.
%
Therefore we can write $X=P\cap C$, where
\bea P=\{(\zeta, Y_2,\ldots, Y_n)\in\Gm^n: f^{\boldsymbol
A}(\zeta^{\power^u}, Y_2^{\power^u},\ldots,Y_n^{\power^u})=0\,,
u=0,\ldots, k-1\} \eea
and
\bea C=\{(\zeta, Y_2,\ldots, Y_n)\in \Gm^n:f^{\boldsymbol
A}(\zeta^{\power^{k+v}},
Y_2^{\power^{k+v}},\ldots,Y_n^{\power^{k+v}})=0\,,
\\v=0,\ldots,T(\power, n,d)-k\} \,.\eea
Now it is clearly enough to show that $C$ is $(\power,
V^{\boldsymbol A})$-stable.

Observe that $\nu=\zeta^{\power^k}$ is a $(\power^{l}-1)$th
root of unity and $\gcd(\power^{l}-1, \power)=1$. Therefore, denoting $W_i=Z(f^{\boldsymbol
A}(\nu^{\power^i},
Y_2^{\power^{k+i}},\ldots,Y_n^{\power^{k+i}}))\subset \Gm^{n-1}$,
we will have
\bea C=\left\{(\zeta, Y_2,\ldots, Y_n)\in \Gm^n: (Y_2,\ldots, Y_n)\in
\bigcap_{i=0}^{l-1} W_i(\power^{l}, T_i^*)\right\}\,\eea
with $T_i^*\ge\lfloor(T(\power, n,
d)-k)/l\rfloor$.

If for each $i$ the subvariety  $S_i=W_i(\power^{l}, T_i^*)$ is $(\power^l,W_i) $-stable, then, clearly, the subvariety
$C$ is $(\power, V^{\boldsymbol A})$-stable.


It is well known (see e. g. Bombieri and Vaaler
\cite{Bombieri-Vaaler}, pp. 27--28) that
$\det(G\cap\Z^{n})=\det(G^\bot\cap\Z^{n})$. Thus $\det(G)=||{\ve
a}||_2$ and (\ref{upper_bound_for_inverse_elements}) together with
(\ref{length-of-a}) implies
\bea S_{f^{\boldsymbol A}}\subset (n\max_{1\le j \le n-1}||{\ve
a}^*_j||_\infty) d B_1^n\varsubsetneq
(nd+n(n-1)\gamma_{n-1}^{\frac{n-1}{2}} d^2)B_1^n\,.\eea
Here $B_1^n$ denotes the unit $n$-ball with respect to the
$l_1$-norm.
Multiplying $f^{\boldsymbol A}$ by a monomial,  we may assume that
$f^{\boldsymbol A}\in \C[Y_1,\ldots,Y_n]$ with 
\bea \deg(f^{\boldsymbol A})< c_2(n,d)=n(n+1)d+\lfloor n(n^2-1)\gamma_{n-1}^{\frac{n-1}{2}} d^2\rfloor\,.\eea
%
%
%
Recall that $k,l\le t_0+ d(d-1)$  by part (ii) of
Lemma \ref{Gcd}. Thus $l T(\power^{l},n-1,
c_2(n,d)\power^{k+l-1})+k\le T(\power,n,d)$ and, consequently, we have $T_i^*\ge T(\power^{l},n-1,
c_2(n,d)\power^{k+l-1})$. Finally, by the inductive
assumption, for each $i$ the subvariety $S_i$ is $(\power^{l},W_i) $-stable.


Suppose now that $U\subset Z(f/g, f^{t_1}/g)\setminus Z(g)$.
Let $r\in \C[X_1,\ldots,X_{n-1}]$ be the resultant of the
polynomials $f/g$ and $f^{t_1}/g$ with respect to the variable
$X_n$. The polynomial $r$ is not identically zero and has the
total degree $\deg(r)\le d^2\power^{t_1}$. The orthogonal
projection $\pi(U)$ into the coordinate subspace corresponding to the indeterminates $X_1,\ldots,X_{n-1}$ clearly lies on the hypersurface $R=Z(r)$
in $\Gm^{n-1}$. In order to apply the inductive assumption we
observe first that $T(\power,n,d)\ge T(\power, n-1,\deg(r))+t_1$
and thus $\pi(U)\subset Q=R(\power, T(\power, n-1,\deg(r)))$.

By the inductive assumption $Q$ is $(\power, R)$-stable and lies
in a finite union $\bigcup D_i$ of $(n-2)$-dimensional torsion
cosets
 in $\Gm^{n-1}$ satisfying conditions of Theorem
\ref{Intersections}. Since $U$ is irreducible, there should be a
torsion coset $D=D_j$ with $\pi(U)\subset D$.
Suppose that $D=Z(h)$ with $h\sim(X_1,\ldots, X_{n-1})^{{\ve a}}- \zeta$. Then the component
$U$  lies in the $n-1$ dimensional coset $D'$ in $\Gm^{n}$ defined by the same polynomial $h$ regarded as a polynomial in $\C[X_1, \ldots, X_n]$.
Moreover the coset $D'$ satisfies conditions of the part (ii) of Theorem \ref{Intersections}. Let us show now that
$U$ is $(\power, V)$-stable.

We have $D={\ve \omega}H_B$, where $B$ is a primitive sublattice
of $\Z^{n-1}$ with $B=\spn_{\Z}({\ve a})$ and $\det(B)=||{\ve
a}||_2\le L(\power, n-1, \deg(r))$.
By Lemma \ref{suitable_matrix}, applied to the subspace
$G=(\spn_{\R}(B))^\bot$, there exists a basis ${\boldsymbol
A}=({\ve a}_1, \ldots, {\ve a}_{n-1})$ of the lattice $\Z^{n-1}$
such that ${\ve a}_1={\ve a}$ and its polar basis ${\boldsymbol
A}^*=({\ve a}^*_1, \ldots, {\ve a}^*_{n-1})$ satisfies the
inequality (\ref{upper_bound_for_inverse_elements}).
%

%
%
%
%
%
The basis ${\boldsymbol A}$ of $\Z^{n-1}$ can be extended to the
basis
\bea {\boldsymbol B}=(({\ve a}_1,0), \ldots, ({\ve a}_{n-1},0),
{\ve e}_n) \eea
of $\Z^{n}$, where $({\ve a}_i,0)$ denotes the vector
$(a_{i1},\ldots,a_{i n-1},0)$ and ${\ve e}_n=(0, \ldots,0, 1)$.

Let $(Y_1,\ldots,Y_n)$ be the coordinates associated with
${\boldsymbol B}$. The rest of the proof follows the scheme used in the previous case with minor changes. However we will give the full proof here for completeness.
As above, we consider the subvariety $X=S^{\boldsymbol B}\cap Z(Y_1-\zeta)$ and observe that
the preimage of $X$ in the initial coordinates is the subvariety
$S\cap Z(h)$. This implies that $U^{\boldsymbol B}$ is an irreducible component
of $X$ and thus it is enough to prove that $X$ is $(\power,
V^{\boldsymbol B})$-stable.

By the inductive assumption $\zeta$ is a
$\power^k(\power^{l}-1)$th root of unity.
%
Therefore, similar to the previous case, we can write $X=P\cap C$, where
\bea P=\{(\zeta, Y_2,\ldots, Y_n)\in\Gm^n: f^{\boldsymbol
B}(\zeta^{\power^u}, Y_2^{\power^u},\ldots,Y_n^{\power^u})=0\,,
u=0,\ldots, k-1\} \eea
and
\bea C=\{(\zeta, Y_2,\ldots, Y_n)\in \Gm^n:f^{\boldsymbol
B}(\zeta^{\power^{k+v}},
Y_2^{\power^{k+v}},\ldots,Y_n^{\power^{k+v}})=0\,,
\\v=0,\ldots,T(\power, n,d)-k\} \,.\eea
The inductive assumption will allow us to show that
 $C$ is
$(\power,V^{\boldsymbol B})$-stable. The stability of $C$ then implies the stability of $X$.

Since $\nu=\zeta^{\power^k}$ is a $(\power^{l}-1)$th
root of unity, it is convenient to introduce the subvarieties  $W_i=Z(f^{\boldsymbol
B}(\nu^{\power^i},
Y_2^{\power^{k+i}},\ldots,Y_n^{\power^{k+i}}))\subset \Gm^{n-1}$ and
to represent $C$ in the form
\bea C=\left\{(\zeta, Y_2,\ldots, Y_n)\in \Gm^n: (Y_2,\ldots, Y_n)\in
\bigcap_{i=0}^{l-1} W_i(\power^{l}, T_i^*)\right\}\,,\eea
where as above  $T_i^*\ge\lfloor(T(\power, n,
d)-k)/l\rfloor$.

Further, the subvariety $C$ is $(\power, V^{\boldsymbol A})$-stable
if each of the subvarieties $S_i=W_i(\power^{l},
T_i^*)$, $0\le i\le l-1$, is
$(\power^l,W_i) $-stable. Thus we need a bound for the size of $\supp(f^{\boldsymbol B})$.

We have $\det(G)=||{\ve a}||_2$ and by
(\ref{upper_bound_for_inverse_elements})
\bea S_{f^{\boldsymbol B}}\subset (n\max_{1\le j \le n-1}||{\ve
a}^*_j||_\infty) d B_1^n\varsubsetneq
(nd+n(n-1)\gamma_{n-1}^{\frac{n-1}{2}} d L(\power, n-1,
\deg(r)))B_1^n\,.\eea
Multiplying $f^{\boldsymbol B}$ by a monomial,  we may assume that
$f^{\boldsymbol B}\in \C[Y_1,\ldots,Y_n]$ with 
\bea \deg(f^{\boldsymbol B})< c_1(\power, n,d)\,.\eea
%
%
%
Recall that $k\le t_0+ d(d-1)$ and $l\le d(d-1)$ by part (ii) of
Lemma \ref{Gcd}. Thus $l T(\power^{l}, n-1, c_1(\power,n,d)\power^{k+l-1})+k\le T(\power, n,d)$ and, by the inductive
assumption, for each $i$ the subvariety $S_i$ is $(\power^{l}, W_i)$-stable.

\section{Proof of Theorem \ref{2D}}

Since the intersection of torsion cosets is a torsion coset itself, we may assume that $V$ is a hypersurface.
The statement of the theorem is clearly true for $n=1$. In this case the subvariety $S(\power, V)$, if not empty, consists of a finite number of roots of unity.

Suppose now that $n=2$ and $V=Z(f)$ with $f\in \C[X_1, X_2]$. Let $S_1$ be an irreducible component of the subvariety $S(\power, V)$.
By the part (ii) of Theorem \ref{Intersections}, the component $S_1$ lies in an one--dimensional torsion coset $D=Z(X_1^{a_1} X_2^{a_2} - \zeta)$ where $\zeta$ is a root of unity, $a_1,a_2\in \Z$
and, by Lemma \ref{lattices_cosets}, we have $\gcd(a_1, a_2)=1$. The integer vector ${\ve a}=(a_1,a_2)$ can be extended to a basis ${\boldsymbol B}=({\ve a}, {\ve b})$
of the lattice $\Z^2$. If $(Y_1, Y_2)$ are the coordinates associated with the basis ${\boldsymbol B}$
then clearly $S_1^{\boldsymbol B}\subset Z(Y_1-\zeta)$. If $Z(Y_1-\zeta)\subset S^{\boldsymbol B}$ then $S_1^{\boldsymbol B}= Z(Y_1-\zeta)$ and the theorem is proved. Thus we may assume without loss of generality that $Z(Y_1-\zeta)\nsubseteq S^{\boldsymbol B}$.

Let $g=f^{\boldsymbol B}$. By the above assumption, for some $u$ the polynomial $Y_1-\zeta$ does not divide the polynomial $g^u$. Next, by part (ii) of Lemma \ref{Gcd}
for some $v$ with $v>u$ all irreducible common factors of the polynomials $g^u$ and $g^v$ define one-dimensional torsion cosets.
Let $h$ be such a factor and suppose that $S_1\subset Z(h)$. Since $h\sim Y_1^{b_1}Y_2^{b_2}-\mu$, where $\mu$ is a root of unity, we will have $S_1=(\zeta,\eta)$ for
some root of unity $\eta$. Thus we only need to settle the case $S_1\nsubseteq Z(\gcd(g^u, g^v))$. Let $r\in \C[Y_2]$ be the resultant of the polynomials $g^u/\gcd(g^u, g^v)$ and $g^v/\gcd(g^u, g^v)$ with respect to the variable $Y_1$. The orthogonal projection of the component $S_1$ into the coordinate axis $Y_2$ clearly lies on the
maximal $(\power, Z(r))$-stable subvariety $S(\power, Z(r))$. Since $S(\power, Z(r))$ is a finite union of roots of unity, the component $S_1$ is a torsion point of $\Gm^2$.

\section{Acknowledgement}

The authors are very grateful to Professor Thomas Tucker for
important comments.

School of Mathematics and Wales Institute of Mathematical and Computational Sciences, Cardiff University, Senghennydd Road, Cardiff CF24 4AG UK

{\em E-mail address:} alievi@cf.ac.uk

School of Mathematics and Maxwell Institute for Mathematical
Sciences, University of Edinburgh, Kings Buildings, Edinburgh EH9
3JZ UK

{\em E-mail address:}  C.Smyth@ed.ac.uk

\enddocument